\documentclass[12pt]{article}
\usepackage{amsmath}
\usepackage{amssymb}
\usepackage{amsfonts}
\usepackage{graphicx}

\begin{document}

\title{\textbf{Semi-cubically hyponormal weighted shifts with Stampfli's
subnormal completion}}
\author{Seunghwan Baek and Mi Ryeong Lee%
\thanks{This research was supported by Basic Science Research
Program through the National Research Foundation of Korea(NRF)
funded by the Ministry of Education (NRF-2017R1D1A1A02017817).}}
\date{}
\maketitle

\begin{abstract} Let $\alpha :1,(1,\sqrt{x},\sqrt{y})^{\wedge }$ be
a weight sequence with Stampfli's subnormal completion and let
$W_{\alpha }$ be its associated weighted shift. In this paper we
discuss some properties of the region
$\mathcal{U}:\mathcal{=}\{(x,y):W_{\alpha }$ is semi-cubically
hyponormal$\}$ and describe the shape of the boundary of
$\mathcal{U}$. In particular, we improve the results of
\cite[Theorem 4.2]{LLB} with properties of $\mathcal{U}$.
\end{abstract}

\bigskip

\centerline{\textbf{1. Introduction}}

\medskip

Let $\mathcal{H}$ be a separable infinite dimensional complex
Hilbert space and let $\mathcal{L}(\mathcal{H})$ be the algebra of
all bounded linear operators on $\mathcal{H}$. A bounded operator
$T$ is said to be \emph{subnormal} if it is the restriction of a
normal operator to an invariant subspace(\cite{Hal}). An operator
$T$ in $\mathcal{L}(\mathcal{H})$ is called \textit{hyponormal} if
$T^{\ast }T\geq TT^{\ast }$. In \cite{Cu0}, Curto defined some
classes of weak subnormality between hyponormality and subnormality in $%
\mathcal{L(H)}$, for examples, $k$-hyponormality and weak
$k$-hyponormality. The weakly $k$-hyponormal weighted shift (whose
definition will be defined below) is the main model in this paper. For a positive integer $k$, an operator $T\in \mathcal{L}(%
\mathcal{H})$ is said to be \textit{weakly} $k$\textit{-hyponormal}
if for every polynomial $p$ of degree $k$ or less, $p(T)$ is
hyponormal (\cite{Cu0},\cite{CuF1},\cite{CuP1},\cite{CuP2}). An
operator $T\in \mathcal{L}(\mathcal{H})$ is called
\textit{semi-weakly} $k$\textit{-hyponormal} if $T+sT^{k}$ is
hyponormal for all complex number $s$(\cite{DEJL}). It is obvious
that a weakly $k$-hyponormal operator is semi-weakly $k$-hyponormal.
In particular,
weak 2-hyponormality is equivalent to semi-weak $2$-hyponormality. The weak $%
2$-hyponormality [weak $3$-hyponormality, semi-weak
$3$-hyponormality, resp.] is referred to as the \textit{quadratic
hyponormality} [\textit{cubic hyponormality,} \textit{semi-cubic
hyponormality, }resp.]. In particular, the quadratic hyponormality
makes an important role in the study of gap on operator properties
such as flatness, completion, and backward extension theory since
1990 (see, for instance,
\cite{BEJ},\cite{Ch},\cite{Cu1},\cite{CuJ},\cite{CuLe},\cite{EJP},\cite{JuP1},\cite{PY}).
In \cite{Cu1}, Curto proved that a $2$-hyponormal weighted shift
with two equal weights $\alpha _{n}=\alpha _{n+1}$ for some
nonnegative integer $n$ has the flatness property, i.e., $\alpha
_{1}=\alpha _{2}=\cdots $. Moreover, he obtained a quadratically
hyponormal weighted shift with first two equal weights which does
not satisfy flatness(\cite{Cu1}). Also in \cite{JuP2}, they showed
that the weighted shift $W_{\alpha }$ with $\alpha
:\sqrt{\frac{2}{3}},\sqrt{\frac{2}{3}},\sqrt{\frac{n+1}{n+2}} $
$(n\geq 2)$ is not cubically hyponormal. Hence the following
question arises naturally(\cite{Cu2}):

\medskip

\textbf{Problem 1.1. }Describe all quadratically hyponormal weighted shifts
with first two equal weights.

\medskip

Recently Li-Cho-Lee in \cite{LCL} proved that if a weighted shift
$W_{\alpha }$ is cubically hyponormal with first two equal weights,
then $W_{\alpha }$ has the flatness property. The structure of
semi-cubically hyponormal weighted shifts has been studied by
several authors(cf.\cite{BEJL1},\cite{BEJL2},\cite{LLB}). To detect
the structure of semi-cubically hyponormal weighted shifts, the
following problem arises naturally:

\medskip

\textbf{Problem 1.2. }Describe all semi-cubically hyponormal weighted shifts
with first two equal weights.

\medskip

As a study of Problem 1.2 it is worthwhile to describe the region $\mathcal{U%
}=\{(x,y):W_{\alpha }$ is semi-cubically hyponormal$\}$ for weighted
shifts $W_{\alpha }$ with weight sequence $\alpha
:1,(1,\sqrt{x},\sqrt{y})^{\wedge },$ where
$(1,\sqrt{x},\sqrt{y})^{\wedge }$ is Stampfli's subnormal
completion. Recall that Curto-Jung studied the shape of the region $
\{(x,y):W_{\alpha }$ is quadratically hyponormal$\}$ in \cite{CuJ}.
In this paper we describe the region $\mathcal{U}$ in detail as a
parallel study.

This note consists of four sections. In Section 2 we recall
characterizations for semi-cubic hyponormality of a weighted shift
$W_{\alpha }$ with weight sequence $\alpha
:1,(1,\sqrt{x},\sqrt{y})^{\wedge }$. In Section 3, we describe the
geometric shapes of the region $\mathcal{U}$ above. In Section 4, we
discuss some remarks concerning the extremality of the region
$\mathcal{U}$.

Throughout this note we denote $\mathbb{R}_{+}$ for the set of
nonnegative real numbers. For a region $\mathcal{V}$ in
$\mathbb{R}^{2}(:= \mathbb{R} \times \mathbb{R})$, we denote the
boundary of $\mathcal{V}$ by $\partial \mathcal{V}$.

Some of the calculations in this paper were aided by using the software tool
\textit{Mathematica} (\cite{Wo}).

\bigskip

\centerline{\textbf{2. Preliminaries}}
\medskip

We recall Stampfli's subnormal completion of three values
(\cite{St}). Let $\alpha _{0},\alpha _{1},\alpha _{2}$ be the first
three weights in $\mathbb{R}_{+}$ such that $\alpha _{0}<\alpha
_{1}<\alpha _{2}$ (to avoid the flatness)\footnote{%
If $W_{\alpha }$ is a subnormal weighted shift such that $\alpha
_{0}=\alpha _{1}$ or $\alpha _{1}=\alpha _{2}$, then $\alpha
_{1}=\alpha _{2}=\cdots .$}. Define
\begin{equation}
\widehat{\alpha }_{n}=\left( \Psi _{1}+\frac{\Psi _{0}}{\alpha _{n-1}^{2}}%
\right) ^{1/2},~\ \ \ \ n\geq 3,  \tag{2.1}
\end{equation}%
where%
\begin{equation*}
\Psi _{0}=-\frac{\alpha _{0}^{2}\alpha _{1}^{2}(\alpha _{2}^{2}-\alpha
_{1}^{2})}{\alpha _{1}^{2}-\alpha _{0}^{2}},\;\;\;\;\;\Psi _{1}=\frac{\alpha
_{1}^{2}(\alpha _{2}^{2}-\alpha _{0}^{2})}{\alpha _{1}^{2}-\alpha _{0}^{2}}.
\end{equation*}%
This produces a bounded sequence $\widehat{\alpha }:=\left\{ \widehat{\alpha
}_{i}\right\} _{i=0}^{\infty },$ where $\widehat{\alpha }_{i}=\alpha _{i}$ $%
(0\leq i\leq 2)$ such that its associated weighted shift
$W_{\widehat{\alpha }}$ is subnormal. As usual, we write $\left(
\alpha _{0},\alpha _{1},\alpha _{2}\right) ^{\wedge }$ for this
weight sequence $\widehat{\alpha }$ induced by (2.1).

We now recall a characterization of the semi-cubic hyponormality of weighted
shifts $W_{\alpha }$ with weight sequence $\alpha :1,(1,\sqrt{x},\sqrt{y}%
)^{\wedge }$.

\medskip

\noindent \textbf{Lemma 2.1 (\cite[Th. 4.1]{LLB}).} \textit{Let }$\alpha :1,(1,\sqrt{x}%
,\sqrt{y})^{\wedge }$\textit{\ with }$1<x<y$\textit{\ be a weight sequence
and let }$W_{\alpha }$\textit{\ be its associated weighted shift. Then }$%
W_{\alpha }$\textit{\ is semi-cubically hyponormal if and only if
}$f\left( x,y\right) :=\sum_{i=0}^{9}\zeta _{i}y^{i}\geq 0,$
\textit{where}
\begin{align*}
\zeta _{0}& =x^{8},~\zeta _{1}=-x^{5}+8x^{6}-18x^{7}+2x^{8}, \\
\zeta _{2}& =x^{2}-8x^{3}+39x^{4}-108x^{5}+131x^{6}-20x^{7}+x^{8}, \\
\zeta _{3}& =-3x+32x^{2}-151x^{3}+338x^{4}-298x^{5}-12x^{6}+10x^{7}, \\
\zeta _{4}& =4-42x+169x^{2}-274x^{3}+40x^{4}+276x^{5}-43x^{6}-4x^{7}, \\
\zeta _{5}& =16x-130x^{2}+359x^{3}-330x^{4}-75x^{5}+34x^{6}, \\
\zeta _{6}& =-2x+38x^{2}-172x^{3}+260x^{4}-34x^{5}-6x^{6}, \\
\zeta _{7}& =-x+4x^{2}+17x^{3}-74x^{4}+18x^{5}, \\
\zeta _{8}& =-2x^{2}+6x^{3}+7x^{4}-2x^{5},~\zeta _{9}=-x^{3}.
\end{align*}

Let $\alpha (x,y):1,(1,\sqrt{x},\sqrt{y})^{\wedge }$ with $1<x<y$ be
a weight sequence with Stampfli's subnormal completion tail and let
$W_{\alpha (x,y)}$ be the
associated weighted shift. For our convenience, we denote $x=1+h$ and $%
y=1+h+k$ $(h,\ k\in \mathbb{R}_{+})$. Then we can rewrite the
polynomials in Lemma 2.1 as following
\begin{equation}
p\left( h,k\right) :=f(1+h,1+h+k)=-\sum_{i=0}^{9}\xi
_{i}(h)k^{i}\geq 0, \tag{2.2}
\end{equation}%
where
\begin{align*}
\xi _{0}(h)& =2h^{9}\left( h+1\right) ^{4}, \
\xi _{1}(h) =h^{8}\left( 16h+7\right) \left( h+1\right) ^{3}, \\
\xi _{2}(h)& =4h^{6}\left( 3h+14h^{2}+14h^{3}-1\right) \left(
h+1\right) ^{2},
\\
\xi _{3}(h)& =h^{5}\left( h+1\right) \left(
3h+98h^{2}+190h^{3}+112h^{4}-4\right) , \\
\xi _{4}(h)& =h^{4}\left( 2h+109h^{2}+322h^{3}+356h^{4}+140h^{5}-5\right) , \\
\xi _{5}(h)& =2h^{3}\left( h+1\right) \left(
5h+46h^{2}+88h^{3}+56h^{4}-1\right) , \\
\xi _{6}(h)& =h^{2}\left( h+1\right) \left(
13h+64h^{2}+104h^{3}+56h^{4}-1\right) , \\
\xi _{7}(h)& =h^{2}\left( h+1\right) \left( 34h+42h^{2}+16h^{3}+9\right) , \\
\xi _{8}(h)& =2h\left( 4h+h^{2}+2\right) \left( h+1\right) ^{2}, \
\xi _{9}(h) =\left( h+1\right) ^{3}.
\end{align*}
For $\alpha (x,y):1,(1,\sqrt{x},\sqrt{y})^{\wedge }$ with $x=1+h$
and $y=1+h+k$ ($h,k \in \mathbb{R}_+$), we denote
\begin{equation*}
\mathcal{R}:=\{(h,k):W_{\alpha (x,y)}\text{ is semi-cubically hyponormal}\}
\end{equation*}%
and%
\begin{equation*}
\mathcal{R}_{\text{q}}:=\{(h,k):W_{\alpha (x,y)}\text{ is quadratically
hyponormal}\}.
\end{equation*}%
 Then it follows from \cite[Theorem 4.2]{LLB} that both $\mathcal{R}\setminus \mathcal{R}_{\text{q}}$
 and $\mathcal{R}_{\text{q}}\setminus \mathcal{R}$ are nonempty sets, \textit{indeed}, a
line segment $\{(\frac{1}{100},k):\beta _{1}\leq k<\alpha _{1}\}$ [or $\{(\frac{%
1}{100},k):\beta _{2}<k\leq \alpha _{2}\}$] contains in $\mathcal{R}\setminus \mathcal{R}_{\text{q}}$
 [or $\mathcal{R}_{\text{q}}\setminus \mathcal{R}$, respectively], where $\alpha _{1}\approx
0.000787776068\cdots $, $\alpha _{2}\approx 0.0422764016\cdots $,
$\beta _{1}\approx 0.000786885627\cdots ,$ and $\beta _{2}\approx
0.0402782805\cdots $; see the proof of \cite[Theorem 4.2]{LLB}. In
the next section, the polynomial in (2.2) can be used to describe the shape of $%
\mathcal{R}$ as a crucial parts.

\bigskip

\centerline{\textbf{3. The shape of the region with semi-cubic
hyponormality}}

\medskip

Let $\alpha (x,y):1,(1,\sqrt{x},\sqrt{y})^{\wedge }$ with $1<x<y$ be
a weight sequence as usual and let $W_{\alpha (x,y)}$ be the
associated weighted shift with $x=1+h$ and $y=1+h+k$ $(h,\ k\in
\mathbb{R}_{+})$. We may replace $k$ by $th$, where $t$ is a
positive real number. Then $p(h,k)$ in (2.2) can be represented by
\begin{equation*}
p(h,k)=p(h,th)=h^{8}\left( \phi _{0}(t)+\phi _{1}(t)h+\phi _{2}(t)h^{2}+\phi
_{3}(t)h^{3}+\phi _{4}(t)h^{4}+\phi _{5}(t)h^{5}\right) ,
\end{equation*}
where
\begin{align*}
 \phi_{0}(t) & =4t^{2}+4t^{3}+5t^{4}+2t^{5}+t^{6}, \\
\phi_{1}(t)&=-2-7t-4t^{2}+t^{3}-2t^{4}-8t^{5}-12t^{6}-9t^{7}-4t^{8}-t^{9}, \\
\phi_{2}(t)&
=-8-37t-76t^{2}-101t^{3}-109t^{4}-102t^{5}-77t^{6}-43t^{7}-16t^{8}-3t^{9}, \\
\phi_{3}(t)&
=-12-69t-180t^{2}-288t^{3}-322t^{4}-268t^{5}-168t^{6}-76t^{7}-22t^{8}-3t^{9},
\\
\phi_{4}(t)&
=-8-55t-168t^{2}-302t^{3}-356t^{4}-288t^{5}-160t^{6}-58t^{7}-12t^{8}-t^{9},
\\
\phi_{5}(t)&
=-2-16t-56t^{2}-112t^{3}-140t^{4}-112t^{5}-56t^{6}-16t^{7}-2t^{8}.
\end{align*}
For brevity, we set
\begin{equation*}
\rho (h,t)=\phi _{0}(t)+\phi _{1}(t)h+\phi _{2}(t)h^{2}+\phi
_{3}(t)h^{3}+\phi _{4}(t)h^{4}+\phi _{5}(t)h^{5}.  \tag{3.1}
\end{equation*}
Then $W_{\alpha(x,y) }$ is semi-cubically hyponormal if and only if
$\rho (h,t)\geq 0$ for all positive numbers $h$ and $t$. We will
detect the set
\begin{equation*}
\mathcal{C}:=\{(h,th)|\rho (h,t)=0\text{ and }h>0,t>0\}\cup \{(0,0)\}
\end{equation*}%
to consider the region of semi-cubic hyponormality of $W_{\alpha
(x,y)}$ below. In fact, the set $\mathcal{C}$ will be a curve (see
Lemma 3.1).

\medskip

\noindent \textbf{Lemma 3.1.}\textit{\ The set
}\emph{$\mathcal{C}$}\textit{\ is a
loop with polar form of }$r=f(\theta )$, $0\leq \theta \leq \frac{\pi }{2}$.%
\textit{\ Therefore }$\mathcal{R}$\textit{\ is a starlike region
with nonempty interior and }$\mathcal{C}=\partial \mathcal{R}$.

\smallskip

\noindent \textit{Proof}. First we fix $t=t_{0}>0$. Since $\phi
_{i}(t_{0}) \ (i=1,...,5)$ are negative obviously,
\begin{equation*}
\frac{\partial }{\partial h}\rho (h,t_{0})=\phi _{1}(t_0)+2\phi
_{2}(t_0)h+3\phi _{3}(t_0)h^{2}+4\phi _{4}(t_0)h^{3}+5\phi
_{5}(t_0)h^{4}
\end{equation*}%
is negative for $h>0$. Then it follows that $\rho (h,t_{0})$ is
decreasing in $h$. Since $\phi _{0}(t_{0})$ is positive, the
equation $\rho (h,t_{0})=0$ of $h$ has a unique solution. $\quad
\blacksquare $

\medskip

The following corollary improves the results of \cite[Theorem
4.2]{LLB}, and its proof follows from the fact that the boundaries
of $\mathcal{R}_{\text{q}}$
and $\mathcal{R}$ are loops with polar forms of\textit{\ }$r=f(\theta )$, $%
0\leq \theta \leq \frac{\pi }{2}$.

\medskip

\noindent \textbf{Corollary 3.2.} \textit{Under the above notation,
we have the following assertions:}

\noindent (i) $\mathcal{R}_{\text{q}}\cap \mathcal{R}$\ \textit{is a
starlike region with nonempty interior,}

\noindent (ii) $\mathcal{R}\setminus \mathcal{R}_{\text{q}}$\textit{\ and }$\mathcal{R}%
_{\text{q}}\setminus \mathcal{R}$\textit{\ are regions with nonempty
interiors.}

\medskip

We now consider the tangent line to the closed curve $\mathcal{C}$ near the
origin.

\medskip

\noindent \textbf{Lemma 3.3.}\textit{\ The tangent line to
}$\mathcal{C}$\textit{\
converges to the }$x$\textit{-axis as }$(h,t)\rightarrow (0^{+},0^{+})$%
\textit{\ and it converges to the }$y$\textit{-axis as }$(k,t)\rightarrow
(0^{+},\infty )$\textit{.}

\smallskip

\noindent \textit{Proof}. We mimic the proof of \cite[Lemma
4.7]{CuJ}. From $k=th$, we have
\begin{equation*}
\frac{dk}{dh}=\frac{dt}{dh} h+t. \tag{3.2} \end{equation*} Since
$\rho (h,t)=0$ on $ \mathcal{C}$, we get
\begin{equation*}
\frac{dt}{dh}=-\frac{\frac{\partial \rho }{\partial
h}}{\frac{\partial \rho }{\partial t}}=-\frac{\phi _{1}(t)+2\phi
_{2}(t)h+3\phi _{3}(t)h^{2}+4\phi _{4}(t)h^{3}+5\phi
_{5}(t)h^{4}}{\phi _{0}^{\prime }(t)+\phi _{1}^{\prime }(t)h+\phi
_{2}^{\prime }(t)h^{2}+\phi _{3}^{\prime }(t)h^{3}+\phi _{4}^{\prime
}(t)h^{4}+\phi _{5}^{\prime }(t)h^{5}}.
\end{equation*}%
Furthermore, we have
\begin{equation*}
\lim_{(h,t)\rightarrow
(0^{+},0^{+})}\frac{dk}{dh}=\lim_{(h,t)\rightarrow
(0^{+},0^{+})}\frac{dt}{dh} h.
\end{equation*}
Since $\lim_{(h,t)\rightarrow (0^{+},0^{+})}\frac{dt}{dh}=\infty $,
by using the L'Hospital's rule and some elementary computations, we
can obtain
\begin{align*}
\lim_{(h,t)\rightarrow (0^{+},0^{+})}\frac{dt}{dh} h&
=\lim_{(h,t)\rightarrow (0^{+},0^{+})}\frac{h}{\frac{dh}{dt}}%
=\lim_{(h,t)\rightarrow (0^{+},0^{+})}\frac{1}{\frac{d}{dh}\left( \frac{dh}{%
dt}\right) } \\
& =\lim_{(h,t)\rightarrow (0^{+},0^{+})}\frac{1}{\frac{\partial }{\partial t}%
\left( \frac{dh}{dt}\right)  \frac{dt}{dh}+\frac{\partial }{\partial h}%
\left( \frac{dh}{dt}\right) } \\
& =\lim_{(h,t)\rightarrow (0^{+},0^{+})}\frac{F_{1}(h,t)}{F_{2}(h,t)}=0,
\end{align*}
for some polynomials $F_{1}$ and $F_{2}$ of $h$ and $t$ (see
Appendix for details) such that
\begin{equation*}
\lim_{(h,t)\rightarrow (0^{+},0^{+})}F_{1}(h,t)=0\text{ and }%
\lim_{(h,t)\rightarrow (0^{+},0^{+})}F_{2}(h,t)=32.
\end{equation*}%
Similarly, we have $\lim_{(k,t)\rightarrow (0^{+},\infty )}\frac{dk}{dh}%
=\infty .$ Hence the proof is complete. $\quad \blacksquare $

\medskip

We now set
\begin{equation*}
h_{M}=\max\{h:(h,k)\in \mathcal{R},~k\in \mathbb{R}_{+}\}, \
k_{M}=\max\{k:(h,k)\in \mathcal{R},~h\in \mathbb{R}_{+}\}. \tag{3.3}
\end{equation*}
Obviously two maximum values $h_{M}$ and $k_{M}$ are well defined.
 Recall that the problem \cite[Problem 5.1]{CuJ} which is finding the values or expressions of
 $\max\{h:(h,k)\in \mathcal{R}_{\text{q}}, \ k\in
\mathbb{R}_{+}\}$ and $\max \{k:(h,k)\in \mathcal{R}_{\text{q}}, \
h\in \mathbb{R}_{+}\}$ are not solved yet. Hence it is worthwhile
finding extremal values $h_{M}$ and $k_{M} $. We discuss the values
of $h_{M}$ and $k_{M}$ below.

\medskip

\noindent \textbf{Lemma 3.4.} \textit{Under the same notation in}
(3.3), \textit{we have} $0<h_{M}<\frac{14}{100}$.

\smallskip

\noindent \textit{Proof.} We can obtain that
\begin{equation*}
\rho \left( \frac{14}{100},t\right) =\frac{1}{156250000}%
\sum_{k=0}^{9}c_{k}t^{k},
\end{equation*}%
where
\begin{eqnarray*}
c_{0} &=&-73892007,c_{1}=-299457081,c_{2}=217020204,c_{3}=195013758, \\
c_{4} &=&243084610,c_{5}=-308008392,c_{6}=-424167096,c_{7}=-364763406, \\
c_{8} &=&-146669607,c_{9}=-32408775.
\end{eqnarray*}
This can be represented by
\begin{align*}
\sum_{k=0}^{9}c_{k}t^{k}& <10^{7}
(-29t+30t^{2}+20t^{3}+25t^{4}-30t^{5}-40t^{6}-35t^{7}-10t^{8}-3t^{9}) \\
& =10^{7}t(-29+30t+20t^{2}+25t^{3}-30t^{4}-40t^{5}-35t^{6}-10t^{7}-3t^{8}) \\
& =10^{7}t((-4+20t^{2}-30t^{4})+(-5+25t^{3}-35t^{6})-10t^{7}-3t^{8}-\eta (t))\\
&=10^{7}t\left(-A^{2}-5B^{2}-\frac{101}{84}-10t^{7}-3t^{8}-\eta(t)\right),
\end{align*}
where $A=\sqrt{30}t^{2}-\sqrt{\frac{10}{3}}$, $B=\sqrt{7}t^{3}-\frac{5}{%
\sqrt{28}}$ and $\eta (t)=20-30t+40t^{5}$. Here, since $\eta (t)$
has exactly one critical number $\sqrt[4]{\frac{3}{20}}$ on
$\mathbb{R_{+}}$ and $\eta ^{\prime \prime }(t)>0$ on
$\mathbb{R_{+}}$, $\eta (t)$ has a positive minimum at
$\sqrt[4]{\frac{3}{20}}$. So, $\rho \left( \frac{14}{100},t\right) $
is negative and since $\mathcal{C}$ is a loop in the first quadrant,
$\mathcal{C}$ lies on the left side of a line $h=\frac{14}{100}$.
$\quad \blacksquare $

\medskip

Recall Descartes' rule of signs that if $p(x)$ is a polynomial with
real coefficients, then the number of positive roots either is equal
to the number of variations in sign of $p(x)$ or is less than that
number by an even number; and the number of negative roots either is
equal to the number of variations in sign of $p(-x)$ or is less than
that number by an even number.

\medskip

\noindent \textbf{Lemma 3.5.} \textit{Given }$h>0$\textit{, there
exist at most two roots} (\textit{possibly a double root})
 $k_{0}>0$ \textit{such that }$p(h,k_{0})=0$\textit{.}

\smallskip

\noindent \textit{Proof}. According to Lemma 3.4, it is sufficient to consider $h<%
\frac{14}{100}$. Recall that
\begin{equation*}
p\left( h,k\right) =-\sum_{i=0}^{9}\xi _{i}(h)k^{i}
\end{equation*}%
where $\xi _{i}(h)$ are shown in (2.2). Here, all coefficients of
$\xi _{0}(h)$, $\xi _{1}(h)$, $\xi _{7}(h)$, $\xi _{8}(h)$, $\xi
_{9}(h)$ are positive and $\xi _{i}(h),i=2,3,4,5,6$ has one
variation in sign, so it has exactly one positive root $\epsilon
_{i},i=2,3,4,5,6$, respectively. Especially, $\epsilon _{6}\approx
0.0584537$ and $\epsilon _{2},\epsilon _{3},\epsilon _{4},\epsilon
_{5}>\frac{14}{100}$, so $\xi _{2}(h),\xi _{3}(h),\xi _{4}(h),\xi
_{5}(h)$ are negative for $h<\frac{14}{100}$. Hence the signs of the
coefficients of $p(h,k)$ change twice as a polynomial in $k$ for
$h<\frac{14}{100}$. By Descartes' rule of signs, it follows that for
fixed $h>0$, the equation $p(h,k)=0$ of $k$ has no or two roots.
$\quad \blacksquare $

\medskip

We may obtain the following lemma similarly.

\medskip

\noindent \textbf{Lemma 3.6.} \textit{Given }$k>0$\textit{, there
exist at most two roots }$($\textit{possibly a double
root}$)$\textit{\ }$h_{0}>0$\textit{\ such that
}$p(h_{0},k)=0$\textit{.}

\medskip

Note that $\mathcal{C}$ consists of two functions $k=f_{1}(h)$ and $%
k=f_{2}(h)$ on the interval $(0,h_{M}]$. Similarly, $\mathcal{C}$
consists of two functions $h=g_{1}(k)$ and $h=g_{2}(k)$ on the
interval $(0,k_{M}]$.

Combining above lemmas, we obtain the main theorem of this paper.

\medskip

\noindent \textbf{Theorem 3.7. }\textit{The region
}$\mathcal{R}$\textit{\ is a simply connected with boundary
}$\partial \mathcal{R}$ \textit{such that}

\noindent (i) $\partial \mathcal{R}$\textit{\ is a loop with polar form }$r=f(\theta )$%
, $0\leq \theta \leq \frac{\pi }{2}$\textit{,}

\noindent (ii) \textit{the tangent lines of} $\partial \mathcal{R}$
\textit{near origin }$(0,0)$\textit{\ converge to the }$x$\textit{-
and }$y$\textit{- axes,}

\noindent (iii) card$(\partial \mathcal{R}\cap \{(a,k):k\in \mathbb{R}\})=2,$ \textit{%
where }$0<a<h_{M}$,

\noindent (iv) card$(\partial \mathcal{R}\cap \{(h,b):h\in \mathbb{R}\})=2,$ \textit{%
where} $0<b<k_{M}$.\footnote{ card$(\cdot)$ denotes for the
cardinality of $\cdot$.}

\bigskip

\centerline{\textbf{4. Further remarks}}

\medskip

Let $\alpha (x,y):1,(1,\sqrt{x},\sqrt{y} )^{\wedge }$ be a weight
sequence with $x=1+h$ and $y=1+h+k$ $(h,\ k\in \mathbb{R}_{+})$. Recall that $h_{M}$ and $%
k_{M}$ is well defined (see Section 3). The problems of expressions
about the extremal values $h_{M}$ and $k_{M}$ are a parallel ones
which were suggested as a question in \cite[Prob. 5.1]{CuJ}. So it
is worth describing the extremal values $h_{M}$ and $k_{M}$. For
this purpose, we denote
\begin{equation*}
Q:=Q(h,t)=\frac{\partial \rho }{\partial t}= \sum_{i=0}^5 \phi
_{i}^{\prime }(t)h^{i}. \tag{4.1}
\end{equation*}
From (3.2), we obtain
\begin{equation*}
\frac{dk}{dh}=\frac{dt}{dh} h + t=\frac{S}{Q},
\end{equation*}%
for a polynomial
\begin{equation}
S:=S(h,t)=\sum_{j=0}^{4}\nu _{j}(t)h^{j},  \tag{4.2}
\end{equation}%
where{\tiny \noindent }%
\begin{eqnarray*}
\nu _{0}(t) &=&8t^{2}+12t^{3}+20t^{4}+10t^{5}+6t^{6}, \\
\nu _{1}(t) &=&2-4t^{2}+2t^{3}-6t^{4}-32t^{5}-60t^{6}-54t^{7}-28t^{8}-8t^{9},
\\
\nu _{2}(t)
&=&16+37t-101t^{3}-218t^{4}-306t^{5}-308t^{6}-215t^{7}-96t^{8}-21t^{9}, \\
\nu _{3}(t)
&=&36+138t+180t^{2}-322t^{4}-536t^{5}-504t^{6}-304t^{7}-110t^{8}-18t^{9}, \\
\nu _{4}(t)
&=&32+165t+336t^{2}+302t^{3}-288t^{5}-320t^{6}-174t^{7}-48t^{8}-5t^{9}, \\
\nu _{5}(t) &=&10+64t+168t^{2}+224t^{3}+140t^{4}-56t^{6}-32t^{7}-6t^{8}.
\end{eqnarray*}
Hence we arrive at the following proposition.

\medskip

\noindent \textbf{Proposition 4.1.} \textit{Under the notation as
in} (3.3), \textit{we have that}

\noindent (i) $h_{M}=\max \left\{ h:\rho (h,t)=0\text{ \textit{and}
}Q(h,t)=0\right\} , $

\noindent (ii) $k_{M}=\max \left\{ th:\rho (h,t)=0\text{ \textit{and} }%
S(h,t)=0\right\} ,$

\noindent \textit{where }$\rho (h,t)$\textit{, }$Q(h,t)$\textit{ and }$%
S(h,t)$ \textit{are as in} (3.1), (4.1) \textit{and }(4.2), \textit{%
respectively. }

\medskip

Before closing this note, we describe the curvature of $\partial
\mathcal{R}$ for the further information above the shape of
$\mathcal{R}$. Since $k=th$,
\begin{equation*}
\frac{d^{2}k}{dh^{2}}=\frac{d^{2}t}{dh^{2}} h+2 \frac{dt}{dh},
\end{equation*}%
and since $\rho (h,t)=0$ on $\mathcal{C}$,
\begin{equation*}
\frac{\partial \rho }{\partial t}\frac{dt}{dh}+\frac{\partial \rho }{%
\partial h}=0.
\end{equation*}%
By differentiation with respect to $h$, we obtain that
\begin{equation*}
\left[ \frac{\partial ^{2}\rho }{\partial t^{2}} \frac{dt}{dh}+\frac{%
\partial }{\partial h}\left( \frac{\partial \rho }{\partial t}\right) \right]
 \frac{dt}{dh}+\frac{\partial \rho }{\partial t} \frac{d^{2}t}{%
dh^{2}}+\frac{\partial }{\partial t}\left( \frac{\partial \rho }{\partial h}%
\right)  \frac{dt}{dh}+\frac{\partial ^{2}\rho }{\partial h^{2}}=0.
\end{equation*}%
Then
\begin{equation*}
\frac{d^{2}t}{dh^{2}}=-\frac{\left( \frac{\partial ^{2}\rho }{\partial t^{2}}%
 \frac{dt}{dh}+\frac{\partial ^{2}\rho }{\partial h\partial t}\right)
 \frac{dt}{dh}+\frac{\partial ^{2}\rho }{\partial t\partial h}
\frac{dt}{dh}+\frac{\partial ^{2}\rho }{\partial
h^{2}}}{\frac{\partial \rho }{\partial t}}. \tag{4.3}
\end{equation*}%
It follows from (4.3) that
\begin{equation*}
\frac{d^{2}k}{dh^{2}} =\frac{2(t+1)P}{Q^{3}},
\end{equation*}
where a polynomial $P:=P(h,t)$ as follows:
\begin{equation*}
P(h,t)=\sum_{j=0}^{14}\mu _{j}(t)h^{j}  \tag{4.4}
\end{equation*}
(see Appendix for detail expression). Hence the curvature $\kappa $ of $%
\mathcal{C}$ can be represented by
\begin{equation*}
\kappa =\frac{\left\vert \frac{d^{2}k}{dh^{2}}\right\vert }{\left( 1+\left(
\frac{dk}{dh}\right) ^{2}\right) ^{\frac{3}{2}}}=\frac{2(t+1)\left\vert
P\right\vert }{\left( Q^{2}+S^{2}\right) ^{\frac{3}{2}}}.
\end{equation*}

\bigskip

\noindent \textbf{Appendix.}

\noindent \textbf{I.} Expressions of polynomials $F_1(h,t)$ and
$F_2(h,t)$ in the proof of Lemma 3.3:

 {\tiny \noindent $F_1(h,t)=32 t + 240 t^2 +
696 t^3 + 1100 t^4 + 1208 t^5 + 1186 t^6 +
 1964 t^7 + 3724 t^8 + 5372 t^9 + 5822 t^{10} + 5200 t^{11} + 4948 t^{12} +
 5760 t^{13} + 6936 t^{14} + 7176 t^{15} + 5976 t^{16} + 3960 t^{17} +
 2066 t^{18} + 836 t^{19} + 252 t^{20} + 52 t^{21} + 6 t^{22}  +
 h (-28 + 312 t + 3949 t^2 + 16915 t^3 + 42682 t^4 + 76734 t^5 +
    110957 t^6 + 144123 t^7 + 182919 t^8 + 230277 t^9 + 280086 t^{10} +
    322338 t^{11} + 346647 t^{12} + 341773 t^{13} + 298967 t^{14} +
    222481 t^{15} + 133414 t^{16} + 58162 t^{17} + 12159 t^{18} - 6127 t^{19} -
    8131 t^{20} - 4849 t^{21} - 1926 t^{22} - 530 t^{23} - 95 t^{24} - 9 t^{25})+
 h^2 (-596 - 2000 t + 12155 t^2 + 106101 t^3 + 387198 t^4 +
    932986 t^5 + 1722473 t^6 + 2654633 t^7 + 3610765 t^8 +
    4460359 t^9 + 5025504 t^{10} + 5104574 t^{11} + 4573489 t^{12} +
    3484739 t^{13} + 2087521 t^{14} + 744375 t^{15} - 224134 t^{16} -
    681226 t^{17} - 708029 t^{18} - 510549 t^{19} - 281857 t^{20} -
    121363 t^{21} - 40200 t^{22} - 9806 t^{23} - 1601 t^{24} - 135 t^{25}) +
 h^3 (-5444 - 45756 t - 152863 t^2 - 192105 t^3 + 342551 t^4 +
    2169529 t^5 + 5602460 t^6 + 10023840 t^7 + 13968978 t^8 +
    15680530 t^9 + 13732205 t^{10} + 7617847 t^{11} - 1748339 t^{12} -
    12002041 t^{13} - 20128503 t^{14} - 23828289 t^{15} - 22604283 t^{16} -
    17834517 t^{17} - 11828954 t^{18} - 6583954 t^{19} - 3043208 t^{20} -
    1145452 t^{21} - 339849 t^{22} - 75179 t^{23} - 11135 t^{24} - 837 t^{25}) +
 h^4 (-28396 - 315428 t - 1641073 t^2 - 5383735 t^3 - 12759905 t^4 -
    24072671 t^5 - 39559920 t^6 - 61420246 t^7 - 94204780 t^8 -
    141860952 t^9 - 203163571 t^{10} - 268956365 t^{11} -
    323508345 t^{12} - 350092895 t^{13} - 338561837 t^{14} -
    290813111 t^{15} - 220432991 t^{16} - 146308305 t^{17} -
    84213036 t^{18} - 41504476 t^{19} - 17215422 t^{20} - 5864066 t^{21} -
    1580995 t^{22} - 317619 t^{23} - 42433 t^{24} - 2835 t^{25}) +
 h^5 (-95064 - 1229024 t - 7622117 t^2 - 30515107 t^3 -
    89806660 t^4 - 210260700 t^5 - 413967832 t^6 - 713119944 t^7 -
    1102750853 t^8 - 1551056775 t^9 - 1992187502 t^{10} -
    2333490762 t^{11} - 2483413660 t^{12} - 2390161932 t^{13} -
    2069286832 t^{14} - 1601598472 t^{15} - 1100170298 t^{16} -
    664814558 t^{17} - 349540997 t^{18} - 157658587 t^{19} -
    59865160 t^{20} - 18639184 t^{21} - 4576000 t^{22} - 831496 t^{23} -
    99425 t^{24} - 5859 t^{25}) +
 h^6 (-217032 - 3115312 t - 21564739 t^2 - 96496565 t^3 -
    316047680 t^4 - 813883800 t^5 - 1728352654 t^6 - 3130041076 t^7 -
    4949119307 t^8 - 6938857825 t^9 - 8704880338 t^{10} -
    9812660254 t^{11} - 9947557956 t^{12} - 9054889260 t^{13} -
    7376411100 t^{14} - 5351445344 t^{15} - 3434783518 t^{16} -
    1933736426 t^{17} - 944299483 t^{18} - 394133661 t^{19} -
    137838348 t^{20} - 39282556 t^{21} - 8754982 t^{22} - 1428540 t^{23} -
    151231 t^{24} - 7749 t^{25}) +
 h^7 (-348936 - 5439080 t - 40904225 t^2 - 198436119 t^3 -
    700907517 t^4 - 1929546419 t^5 - 4327281625 t^6 -
    8154156563 t^7 - 13200588777 t^8 - 18650312475 t^9 -
    23240249066 t^{10} - 25705008734 t^{11} - 25310883394 t^{12} -
    22192965638 t^{13} - 17292661250 t^{14} - 11925506814 t^{15} -
    7234187498 t^{16} - 3827432934 t^{17} - 1746078261 t^{18} -
    676372203 t^{19} - 217855521 t^{20} - 56650759 t^{21} - 11386853 t^{22} -
    1650831 t^{23} - 152261 t^{24} - 6615 t^{25})  +
 h^8 (-400344 - 6694936 t - 53947543 t^2 - 279612425 t^3 -
    1050175729 t^4 - 3053331091 t^5 - 7168770397 t^6 -
    13995538595 t^7 - 23204634559 t^8 - 33178326413 t^9 -
    41355871918 t^{10} - 45261502546 t^{11} - 43668664482 t^{12} -
    37188395942 t^{13} - 27919995234 t^{14} - 18414048926 t^{15} -
    10605306830 t^{16} - 5287921306 t^{17} - 2255578187 t^{18} -
    809818005 t^{19} - 239315517 t^{20} - 56403895 t^{21} - 10119601 t^{22} -
    1283279 t^{23} - 100555 t^{24} - 3537 t^{25}) +
 h^9 (-326796 - 5819448 t - 49837988 t^2 - 273713140 t^3 -
    1084717980 t^4 - 3309293116 t^5 - 8097011048 t^6 -
    16340639312 t^7 - 27752587844 t^8 - 40252817208 t^9 -
    50388715192 t^{10} - 54833893976 t^{11} - 52097950568 t^{12} -
    43289239112 t^{13} - 31431822784 t^{14} - 19876680400 t^{15} -
    10881899300 t^{16} - 5111484152 t^{17} - 2034098052 t^{18} -
    673895764 t^{19} - 181422108 t^{20} - 38344924 t^{21} - 6044408 t^{22} -
    654144 t^{23} - 41708 t^{24} - 1080 t^{25}) +
  h^{10} (-185764 - 3504736 t - 31732572 t^2 - 183720188 t^3 -
    764643220 t^4 - 2438414988 t^5 - 6200763080 t^6 -
    12918995824 t^7 - 22480185500 t^8 - 33127722240 t^9 -
    41757346984 t^{10} - 45330176584 t^{11} - 42553774680 t^{12} -
    34599932712 t^{13} - 24345257600 t^{14} - 14772424304 t^{15} -
    7681584172 t^{16} - 3390268976 t^{17} - 1252661276 t^{18} -
    380104860 t^{19} - 92195188 t^{20} - 17189932 t^{21} - 2321048 t^{22} -
    205376 t^{23} - 9780 t^{24} - 144 t^{25}) +
 h^{11} (-70004 - 1394044 t - 13295268 t^2 - 80865340 t^3 -
    352387812 t^4 - 1171762044 t^5 - 3091771444 t^6 -
    6645314508 t^7 - 11850492600 t^8 - 17764157752 t^9 -
    22590893768 t^{10} - 24522310008 t^{11} - 22800228776 t^{12} -
    18176491800 t^{13} - 12406481544 t^{14} - 7220796920 t^{15} -
    3558313860 t^{16} - 1468817100 t^{17} - 500096660 t^{18} -
    137397580 t^{19} - 29512692 t^{20} - 4726252 t^{21} - 522564 t^{22} -
    34620 t^{23} - 976 t^{24})+
 h^{12} (-15740 - 329860 t - 3304580 t^2 - 21062780 t^3 -
    95900500 t^4 - 331982060 t^5 - 907976940 t^6 - 2012593620 t^7 -
    3679370520 t^8 - 5616001000 t^9 - 7216427240 t^{10} -
    7847259160 t^{11} - 7239579880 t^{12} - 5666569880 t^{13} -
    3753583000 t^{14} - 2093078760 t^{15} - 974096940 t^{16} -
    373547220 t^{17} - 115876820 t^{18} - 28307500 t^{19} - 5230340 t^{20} -
    684860 t^{21} - 56380 t^{22} - 2180 t^{23}) +
 h^{13} (-1600 - 35200 t - 369600 t^2 - 2464000 t^3 - 11704000 t^4 -
    42134400 t^5 - 119380800 t^6 - 272870400 t^7 - 511632000 t^8 -
    795872000 t^9 - 1034633600 t^{10} - 1128691200 t^{11} -
    1034633600 t^{12} - 795872000 t^{13} - 511632000 t^{14} -
    272870400 t^{15} - 119380800 t^{16} - 42134400 t^{17} - 11704000 t^{18} -
    2464000 t^{19} - 369600 t^{20} - 35200 t^{21} - 1600 t^{22}) $}

    \smallskip

{\tiny \noindent $F_2(h,t)=32 + 64 t - 1032 t^2 - 6328 t^3 - 23636
t^4 - 62984 t^5 -
 126878 t^6 - 204834 t^7 - 271766 t^8 - 305814 t^9 - 299656 t^{10} -
 264344 t^{11} - 216716 t^{12} - 168156 t^{13} - 122788 t^{14} - 81812 t^{15} -
 47956 t^{16} - 23728 t^{17} - 9542 t^{18} - 2946 t^{19} - 638 t^{20} -
 78 t^{21} +
 h (676 + 3586 t + 2334 t^2 - 29494 t^3 - 153650 t^4 - 457552 t^5 -
    977120 t^6 - 1620248 t^7 - 2177740 t^8 - 2444682 t^9 -
    2361878 t^{10} - 2026450 t^{11} - 1585658 t^{12} - 1134518 t^{13} -
    709066 t^{14} - 339474 t^{15} - 69254 t^{16} + 74324 t^{17} +
    108236 t^{18} + 82180 t^{19} + 43808 t^{20} + 17178 t^{21} + 4854 t^{22} +
    910 t^{23} + 90 t^{24}) +
 h^2 (5972 + 43324 t + 131516 t^2 + 203974 t^3 + 42976 t^4 -
    671308 t^5 - 2084212 t^6 - 3834008 t^7 - 5102506 t^8 -
    5135662 t^9 - 3790318 t^{10} - 1569368 t^{11} + 838184 t^{12} +
    2948152 t^{13} + 4478564 t^{14} + 5202270 t^{15} + 4992836 t^{16} +
    4002220 t^{17} + 2665056 t^{18} + 1456828 t^{19} + 641762 t^{20} +
    221154 t^{21} + 56618 t^{22} + 9736 t^{23} + 864 t^{24}) +
 h^3 (31128 + 276836 t + 1150780 t^2 + 3038924 t^3 + 5842372 t^4 +
    8936220 t^5 + 11926996 t^6 + 15668872 t^7 + 22209040 t^8 +
    33176468 t^9 + 47998028 t^{10} + 63676968 t^{11} + 76223768 t^{12} +
    82328528 t^{13} + 80261184 t^{14} + 70230164 t^{15} + 54519516 t^{16} +
    36970752 t^{17} + 21512160 t^{18} + 10521948 t^{19} + 4216244 t^{20} +
    1334988 t^{21} + 315284 t^{22} + 49828 t^{23} + 3996 t^{24}) +
 h^4 (111376 + 1163964 t + 5821476 t^2 + 18852526 t^3 +
    45270620 t^4 + 87779216 t^5 + 146456144 t^6 + 220847572 t^7 +
    310821406 t^8 + 413240770 t^9 + 516947748 t^{10} + 601932518 t^{11} +
    645387458 t^{12} + 631503014 t^{13} + 559292366 t^{14} +
    444297688 t^{15} + 313149838 t^{16} + 193256978 t^{17} +
    102763740 t^{18} + 46145440 t^{19} + 17038490 t^{20} + 4979066 t^{21} +
    1083454 t^{22} + 156704 t^{23} + 11340 t^{24}) +
 h^5 (293464 + 3500016 t + 20100912 t^2 + 74817136 t^3 +
    205306596 t^4 + 448518102 t^5 + 822643010 t^6 + 1316777034 t^7 +
    1889956338 t^8 + 2472431338 t^9 + 2968407886 t^{10} +
    3270564242 t^{11} + 3293021370 t^{12} + 3010938526 t^{13} +
    2481416882 t^{14} + 1827552430 t^{15} + 1190843358 t^{16} +
    678229278 t^{17} + 332524986 t^{18} + 137590838 t^{19} +
    46762762 t^{20} + 12547860 t^{21} + 2495060 t^{22} + 326992 t^{23} +
    21168 t^{24}) +
 h^6 (583256 + 7750260 t + 49660748 t^2 + 205752514 t^3 +
    624207304 t^4 + 1489683736 t^5 + 2933740520 t^6 +
    4933506878 t^7 + 7260575274 t^8 + 9509362730 t^9 +
    11196187778 t^{10} + 11903522618 t^{11} + 11429636038 t^{12} +
    9881660930 t^{13} + 7651486738 t^{14} + 5267652626 t^{15} +
    3194619850 t^{16} + 1686802730 t^{17} + 763835850 t^{18} +
    290747772 t^{19} + 90475018 t^{20} + 22093662 t^{21} + 3965742 t^{22} +
    464120 t^{23} + 26460 t^{24})+
 h^7 (873856 + 12700836 t + 88969692 t^2 + 401747196 t^3 +
    1320338460 t^4 + 3381881092 t^5 + 7059190828 t^6 +
    12391350324 t^7 + 18712526596 t^8 + 24709251912 t^9 +
    28843701320 t^{10} + 29956295656 t^{11} + 27750750696 t^{12} +
    22913898632 t^{13} + 16805073656 t^{14} + 10881380680 t^{15} +
    6167473560 t^{16} + 3024886420 t^{17} + 1264241468 t^{18} +
    441023228 t^{19} + 124740796 t^{20} + 27412596 t^{21} + 4373756 t^{22} +
    448036 t^{23} + 21924 t^{24}) +
 h^8 (974768 + 15295716 t + 115500420 t^2 + 560373234 t^3 +
    1968384896 t^4 + 5348978548 t^5 + 11733230762 t^6 +
    21397707568 t^7 + 33141761300 t^8 + 44275237192 t^9 +
    51575870680 t^{10} + 52756294052 t^{11} + 47556062480 t^{12} +
    37797612696 t^{13} + 26425566852 t^{14} + 16167339248 t^{15} +
    8585647192 t^{16} + 3912087476 t^{17} + 1505372468 t^{18} +
    478610362 t^{19} + 121913648 t^{20} + 23778836 t^{21} + 3305794 t^{22} +
    287888 t^{23} + 11556 t^{24}) +
 h^9 (792996 + 13315534 t + 107368306 t^2 + 554427366 t^3 +
    2063172726 t^4 + 5903378738 t^5 + 13531539302 t^6 +
    25555118362 t^7 + 40571735190 t^8 + 54943226572 t^9 +
    64128413172 t^{10} + 64956009820 t^{11} + 57314726300 t^{12} +
    44094205476 t^{13} + 29520345996 t^{14} + 17114452244 t^{15} +
    8522622640 t^{16} + 3602070918 t^{17} + 1270507066 t^{18} +
    365266910 t^{19} + 82772814 t^{20} + 14067338 t^{21} + 1656174 t^{22} +
    116914 t^{23} + 3510 t^{24}) +
 h^{10} (455412 + 8132992 t + 69585864 t^2 + 380066784 t^3 +
    1489743296 t^4 + 4466438104 t^5 + 10659299952 t^6 +
    20803140872 t^7 + 33840253844 t^8 + 46512317016 t^9 +
    54538481400 t^{10} + 54900887704 t^{11} + 47607247272 t^{12} +
    35584281768 t^{13} + 22876433512 t^{14} + 12583424136 t^{15} +
    5871001220 t^{16} + 2293517384 t^{17} + 736427040 t^{18} +
    189310760 t^{19} + 37503928 t^{20} + 5402976 t^{21} + 514120 t^{22} +
    26832 t^{23} + 468 t^{24}) +
 h^{11} (174600 + 3301624 t + 29843624 t^2 + 171694488 t^3 +
    706230840 t^4 + 2211748200 t^5 + 5483321048 t^6 +
    11045118280 t^7 + 18406293584 t^8 + 25700549680 t^9 +
    30329986960 t^{10} + 30417361520 t^{11} + 25990346096 t^{12} +
    18917805904 t^{13} + 11694584240 t^{14} + 6101826320 t^{15} +
    2660458792 t^{16} + 955087448 t^{17} + 276334408 t^{18} +
    62465400 t^{19} + 10527960 t^{20} + 1226696 t^{21} + 86008 t^{22} +
    2600 t^{23}) +
 h^{12} (40064 + 799472 t + 7610256 t^2 + 45985712 t^3 +
    198009488 t^4 + 646511328 t^5 + 1662866016 t^6 + 3454973376 t^7 +
    5899023936 t^8 + 8374216864 t^9 + 9960028000 t^{10} +
    9968358816 t^{11} + 8407271392 t^{12} + 5966724736 t^{13} +
    3547563264 t^{14} + 1753006272 t^{15} + 711137472 t^{16} +
    232573680 t^{17} + 59706512 t^{18} + 11555888 t^{19} + 1579536 t^{20} +
    135200 t^{21} + 5408 t^{22}) +
 h^{13} (4160 + 87360 t + 873600 t^2 + 5532800 t^3 + 24897600 t^4 +
    84651840 t^5 + 225738240 t^6 + 483724800 t^7 + 846518400 t^8 +
    1222748800 t^9 + 1467298560 t^{10} + 1467298560 t^{11} +
    1222748800 t^{12} + 846518400 t^{13} + 483724800 t^{14} +
    225738240 t^{15} + 84651840 t^{16} + 24897600 t^{17} + 5532800 t^{18} +
    873600 t^{19} + 87360 t^{20} + 4160 t^{21}) $}

\medskip

\noindent \textbf{II.} Expressions of polynomials in (4.4):

{\tiny \noindent $\mu
_{0}(t)=128t^{2}+704t^{3}+1824t^{4}+3408t^{5}+4336t^{6}+4896t^{7}+5128t^{8}+6804t^{9}+9068t^{10}+10988t^{11}+10708t^{12}+8392t^{13}+5200t^{14}+2500t^{15}+916t^{16}+228t^{17}+36t^{18}
$ }

{\tiny \noindent $\mu
_{1}(t)=-16-256t+388t^{2}+7732t^{3}+33538t^{4}+89312t^{5}+173263t^{6}+267617t^{7}+346987t^{8}+388575t^{9}+378876t^{10}+315204t^{11}+213982t^{12}+106838t^{13}+25618t^{14}-15502t^{15}-23942t^{16}-16828t^{17}-7933t^{18}-2623t^{19}-569t^{20}-69t^{21}
$ }

{\tiny \noindent $\mu
_{2}(t)=-240-4300t-19140t^{2}-37711t^{3}-9697t^{4}+146455t^{5}+450429t^{6}+772669t^{7}+878315t^{8}+587441t^{9}-83457t^{10}-911958t^{11}-1584574t^{12}-1858042t^{13}-1683170t^{14}-1213193t^{15}-689391t^{16}-294345t^{17}-80655t^{18}-2459t^{19}+10703t^{20}+6183t^{21}+1945t^{22}+364t^{23}+36t^{24}
$ }

{\tiny \noindent $\mu
_{3}(t)=-1166-24710t-163886t^{2}-616343t^{3}-1602434t^{4}-3232844t^{5}-5570122t^{6}-8832384t^{7}-13254901t^{8}-18479035t^{9}-23128813t^{10}-25285116t^{11}-23670276t^{12}-18574828t^{13}-11764294t^{14}-5461319t^{15}-1181228t^{16}+801896t^{17}+1160836t^{18}+814598t^{19}+397517t^{20}+142173t^{21}+36499t^{22}+6184t^{23}+540t^{24}
$ }

{\tiny \noindent $\mu
_{4}(t)=-842-50304t-468240t^{2}-2265997t^{3}-7308965t^{4}-17630785t^{5}-34118895t^{6}-55457271t^{7}-77631099t^{8}-94005675t^{9}-96947097t^{10}-81696948t^{11}-50483568t^{12}-12977216t^{13}+18418352t^{14}+35271673t^{15}+36834749t^{16}+28571797t^{17}+17585255t^{18}+8719399t^{19}+3457271t^{20}+1068063t^{21}+244121t^{22}+37168t^{23}+2862t^{24}
$ }

{\tiny \noindent $\mu
_{5}(t)=13640+92394t+137046t^{2}-868945t^{3}-5619738t^{4}-17107840t^{5}-33978638t^{6}-45733236t^{7}-33259259t^{8}+22915663t^{9}+129641272t^{10}+270547173t^{11}+406486747t^{12}+491304561t^{13}+495874973t^{14}+424072670t^{15}+308635149t^{16}+190776423t^{17}+99416612t^{18}+43090742t^{19}+15199531t^{20}+4213351t^{21}+865207t^{22}+117524t^{23}+7938t^{24}
$ }

{\tiny \noindent $\mu
_{6}(t)=66520+810798t+4814770t^{2}+18907409t^{3}+56535319t^{4}+140438049t^{5}+305526531t^{6}+595876284t^{7}+1045146056t^{8}+1640428858t^{9}+2291862442t^{10}+2840001732t^{11}+3113451244t^{12}+3011554912t^{13}+2560812332t^{14}+1904502598t^{15}+1230170626t^{16}+683642376t^{17}+322702188t^{18}+127102163t^{19}+40716845t^{20}+10206127t^{21}+1878921t^{22}+225718t^{23}+13230t^{24}
$ }

{\tiny \noindent $\mu
_{7}(t)=162236+2312142t+16067730t^{2}+73104575t^{3}+246765666t^{4}+663242262t^{5}+1483787440t^{6}+2841803261t^{7}+4739395375t^{8}+6948282375t^{9}+8997012383t^{10}+10307245651t^{11}+10444972061t^{12}+9345051607t^{13}+7355999639t^{14}+5067050623t^{15}+3031227839t^{16}+1558464919t^{17}+678897523t^{18}+245692946t^{19}+71828821t^{20}+16268871t^{21}+2668789t^{22}+280256t^{23}+13986t^{24}
$ }

{\tiny \noindent $\mu
_{8}(t)=251284+3925470t+29755794t^{2}+146383587t^{3}+527410783t^{4}+1488333325t^{5}+3432113331t^{6}+6652948356t^{7}+11048757104t^{8}+15917079254t^{9}+20043961474t^{10}+22153678856t^{11}+21519567732t^{12}+18354629988t^{13}+13704662132t^{14}+8911130946t^{15}+5006861682t^{16}+2404435724t^{17}+972040780t^{18}+323881861t^{19}+86296013t^{20}+17575887t^{21}+2545529t^{22}+229834t^{23}+9450t^{24}
$ }

{\tiny \noindent $\mu
_{9}(t)=264248+4424166t+35788630t^{2}+186720973t^{3}+707640008t^{4}+2079733506t^{5}+4939315353t^{6}+9746624516t^{7}+16290357438t^{8}+23369615740t^{9}+29027484452t^{10}+31379283482t^{11}+29587502248t^{12}+24324403820t^{13}+17387207738t^{14}+10749106608t^{15}+5700898724t^{16}+2563914798t^{17}+962079758t^{18}+294436673t^{19}+71134864t^{20}+12919554t^{21}+1629861t^{22}+123428t^{23}+3942t^{24}
$ }

{\tiny \noindent $\mu
_{10}(t)=191880+3404534t+29082298t^{2}+159462286t^{3}+631314510t^{4}+1924393206t^{5}+4701868130t^{6}+9461239838t^{7}+15978840550t^{8}+22951973380t^{9}+28292750396t^{10}+30096283108t^{11}+27697792564t^{12}+22049931956t^{13}+15142201948t^{14}+8920331444t^{15}+4469046812t^{16}+1880437606t^{17}+652963306t^{18}+182547246t^{19}+39649374t^{20}+6338550t^{21}+682082t^{22}+41614t^{23}+918t^{24}
$ }

{\tiny \noindent $\mu
_{11}(t)=94946+1772384t+15880692t^{2}+90982752t^{3}+374575704t^{4}+1180659948t^{5}+2963589944t^{6}+6082382852t^{7}+10395813954t^{8}+14988362612t^{9}+18389054228t^{10}+19303041588t^{11}+17379547484t^{12}+13418145020t^{13}+8856954252t^{14}+4968428396t^{15}+2346437330t^{16}+920274780t^{17}+293998536t^{18}+74430684t^{19}+14343420t^{20}+1976552t^{21}+174956t^{22}+7968t^{23}+90t^{24}
$ }

{\tiny \noindent $\mu
_{12}(t)=30566+596842t+5579030t^{2}+33236026t^{3}+141717410t^{4}+460464766t^{5}+1184963698t^{6}+2477876174t^{7}+4285271260t^{8}+6204406020t^{9}+7582152252t^{10}+7859337252t^{11}+6924135972t^{12}+5181171804t^{13}+3281281860t^{14}+1747111420t^{15}+774002078t^{16}+281000242t^{17}+81795982t^{18}+18493250t^{19}+3094618t^{20}+354182t^{21}+23914t^{22}+662t^{23}
$ }

{\tiny \noindent $\mu
_{13}(t)=5760+117000t+1134840t^{2}+6993960t^{3}+30741240t^{4}+102531600t^{5}+269537040t^{6}+572562720t^{7}+999542880t^{8}+1450476720t^{9}+1762572240t^{10}+1800867120t^{11}+1548927120t^{12}+1119698880t^{13}+677214720t^{14}+340002720t^{15}+140001120t^{16}+46450440t^{17}+12099960t^{18}+2379240t^{19}+331320t^{20}+29040t^{21}+1200t^{22}
$ }

{\tiny \noindent $\mu
_{14}(t)=480+10080t+100800t^{2}+638400t^{3}+2872800t^{4}+9767520t^{5}+26046720t^{6}+55814400t^{7}+97675200t^{8}+141086400t^{9}+169303680t^{10}+169303680t^{11}+141086400t^{12}+97675200t^{13}+55814400t^{14}+26046720t^{15}+9767520t^{16}+2872800t^{17}+638400t^{18}+100800t^{19}+10080t^{20}+480t^{21}
$ }

\bigskip

\bigskip

Seunghwan Baek \vskip0.1cm  Faculty of Liberal Education,

Kyungpook National University,

Daegu 41566, Korea

\textit{E-mail address}: \textit{seunghwan@knu.ac.kr.}

\bigskip

Mi Ryeong Lee \vskip 0.1cm Institute of Liberal Education,

Daegu Catholic University,

Gyeongsan, Gyeongbuk 38430, Korea

\textit{E-mail address}: \textit{leemr@cu.ac.kr.}

\end{document}